\documentclass[12pt,leqno]{article}

\usepackage{amsmath} 
\usepackage{amssymb}
\usepackage[english]{babel}
\usepackage{latexsym}

\newfont{\graf}{eufm10}

\newcommand{\be}{\begin{equation}}
\newcommand{\ee}{\end{equation}}
\newcommand{\bea}{\begin{eqnarray}}
\newcommand{\eea}{\end{eqnarray}}
\newcommand{\bean}{\begin{eqnarray*}}
\newcommand{\eean}{\end{eqnarray*}}

\newcommand{\CD}{{\cal D}}
\newcommand{\TT}{\mathbb{T}}
\newcommand{\CC}{\mathbb{C}}
\newcommand{\DD}{\mathbb{D}}
\newcommand{\BB}{\mathbb{B}}

\newcommand{\ep}{\varepsilon}

\newcounter{alteqn}

\title{A Cantor set in the unit sphere in $\mathbb{C}^2$ with large
polynomial hull \footnote{
sets, polynomial and rational hulls, removable singularities of CR
functions.}
\footnote{{\it email} joericke@math.uu.se}}

\author{Burglind J\"oricke}

\begin{document}

\maketitle
\thispagestyle{plain}

\begin{abstract} It is an old question how massive polynomial
hulls of Cantor sets in $\mathbb{C}^n$ can be. In contrast to
expectation e.g. Rudin, Vitushkin and Henkin showed on examples that
it can be rather massive. Motivated by problems of holomorphic convexity
of subsets of strictly pseudoconvex boundaries and removable
singularities the question was asked for Cantor sets in the unit
sphere. It was known that tame Cantor sets in the unit sphere are
polynomially convex. We give an example of a wild Cantor set in the
sphere whose polynomial hull contains a large ball. In some sense this can
be opposed to a still open conjecture of Vitushkin on the existence of a
lower bound for the diameter of the largest boundary component of a
relatively closed complex curve in the ball passing through the origin.

\end{abstract}


It is an old question by Walter Rudin asked in connection with Banach algebras
and approximation by polynomials how massive polynomial hulls of Cantor sets
may be. In contrast to expectation (note that Cantor sets have Hausdorff
dimension zero) e.g. Rudin, Vitushkin and Henkin showed on examples that the
mentioned polynomial hull can be rather massive. Rudin himself constructed a
Cantor set in $\mathbb C ^2$ such that its polynomial hull (and even its
rational hull) contains an analytic variety of dimension one (\cite{R},
theorem 5; \cite{G}, theorem III.2.5 ). Later Vitushkin \cite{V} and Henkin
\cite{H} gave examples of Cantor sets with interior points in the polynomial
hull. 

The problem received new attention in connection with interest in topology
on strictly pseudoconvex boundaries and hulls of their subsets as well as in
connection with removable singularities of CR functions. In particular, it was
asked whether Cantor sets in the unit sphere $\mathbb C ^2$ are polynomially
convex. The expectation was that for subsets of the sphere the situation
changes dramatically as it is the case for some other problems. E.g. totally
real discs in $\mathbb C ^2$ are not necessarily polynomially convex
\cite{NW}, but if contained in the sphere they are so \cite{J}. Further, the
polynomial hull of a compact set in $\mathbb C ^2$ of finite one-dimensional
Hausdorff measure is not necessarily an analytic variety \cite{A}, but if the
set is contained in the sphere this is so \cite{A}. Moreover, the question
about polynomial hulls of Cantor sets in the sphere has some relation to a
still open conjecture of Vitushkin on the existence of a lower bound for the
diameter of the largest boundary component of a relatively closed complex
curve in the ball passing through the origin.  

In \cite{FS} the slightly more general question was raised whether Cantor sets
in boundaries $\partial G$ of strictly pseudoconvex domains $G$ in $\mathbb C
^2$ are convex with respect to the space of holomorphic functions in $G$ which
are continuous in the closure $\overline G$ of the domain $G$.  The question
was answered in \cite{FS} affirmatively for some class of Cantor sets. The
main tools used for this in \cite{FS} are the theorem of Bedford and
Klingenberg \cite{BK}, a characterization of tame Cantor sets in $\mathbb R
^n$ due to Bing and the non-trivial fact that a $C^2$ manifold that is
homeomorphic to $\mathbb R ^3$ is $C^2$ diffeomorphic to $\mathbb R ^3$ (see
references in \cite{FS}). With the mentioned
tools in mind it was a natural step to obtain the following result.\\


\noindent {\bf Theorem a.} {\it Tame Cantor sets in the unit sphere in $\mathbb
 C ^2$ are polynomially convex.} \\

The theorem follows directly from \cite{E} and the mentioned result of
Bing and was proved independently a bit later by M. Lawrence \cite{L}.

Recall that a compact set $E \subset \mathbb{R}^n$ is a tame Cantor
set if there is a homeomorphism of $\mathbb{R}^n$ which carries $E$ to
the middle-third Cantor set in a coordinate line. Bing's result is
that this is equivalent to the following separation property.

For each pair of different points $p$ and $q$  in $E$ and any $\ep >0$
there exists a set $\overline{b^n} \subset\mathbb R ^n$ homeomorphic
to a closed n-ball of diameter less than $\ep$ with boundary $\partial
b^n = {\cal S} ^{n-1}$ disjoint from $E$ with $p\in b^n$ and $q\notin
\overline{b^n}$.

Theorem a can be stated in a slightly more general form. For a
domain $G \subset \mathbb C ^n$ denote by $A(G)$ the algebra of
analytic functions in $G$ which are continuous in $\overline G$.  For
a compact subset $K$ of $\overline G$ we consider its $A(G)$-hull,
\[
A(G)-\mbox{hull}(K) = \{z \in \overline G : |f(z)|\le \max_{K} |f|
\mbox{ for all functions}\;  f \in A(G) \}.
\]

In case $G = \mathbb C ^n$ we obtain the polynomial hull
\[
\hat K = \{z \in \mathbb C ^n : |p(z)| \le \max_{K} |p|  \mbox{ for
  all polynomials } p \}
\]
of the compact subset $K$ of  $\mathbb C^n$. 

The rational hull of the compact $K$ is defined by replacing in the
last condition polynomials by rational functions which are analytic in
a neighbourhood of the compact $K$.\\

\noindent {\bf Theorem b.} {\it Let $G$ be a bounded strictly pseudoconvex
 domain in $\mathbb C^2$ with smooth boundary. Let $K \subset
 \partial G$ be a Cantor set with the following separation property.

For each pair of distinct points $p$ and $q$  in $K$ there exists a
smooth 2-sphere ${\cal S}^2 \subset \partial G \setminus K$ which
separates $p$ and $q$ (i.e. $p$ and $q$ are in different connected
components of $\partial G \setminus {\cal S}^2 $).

Then $K$ is $A(G)$-convex (i.e. $A(G)-\mbox{hull}(K)=K$).}\\

Theorem a follows from theorem b by the result of Bing and the
mentioned fact about homeomorphic and diffeomorphic $\mathbb R ^3$
(for more detail see \cite{FS}). Note that it is not necessary to
require in theorem b that the 2-spheres bound 3-balls in $\partial G$
and that these balls have small size. In the description of tame
Cantor sets the condition that the $b ^3$ has small diameter is
important. There are wild Cantor sets in $\mathbb R ^3$ (\cite{Bi})
each pair of distinct points of which can be separated by spheres (but
not by spheres of small size). 

The proof of theorem b follows from the theorem of Bedford and
Klingenberg. Here is a sketch of its proof for convenience of the
reader (see \cite{FS}, \cite{E}, \cite{L}).

If $K$ were not $A(G)$-convex, then by Zorn's lemma there would be a
minimal compact subset $K'$ of $K$ containing a given point of $G$ in
its $A(G)$-hull. Then $A(G)$-hull$(K')$ is connected by Rossi's
maximum principle. By the separation property there is a smooth
2-sphere ${\cal S}^2 \subset \partial G \setminus K'$ which divides
$K'$ and hence also $\partial G$. By \cite{BK} (maybe after perturbing
the 2-sphere) ${\cal S}^2$ bounds a Levi-flat 3-ball ${\cal B}^3
\subset G$ such that $\overline {\cal B}^3$ is equal to the envelope
of holomorphy of ${\cal S}^2$ and to $A(G)$-hull$({\cal S}^2)$. By a
theorem of Alexander and Stout (see the references in \cite{St})
$\overline {\cal B}^3$ divides $\overline G$. The envelope of
holomorphy of $\partial G \setminus K'$ equals $\overline G \setminus
A(G)- \mbox{hull}(K')$ (see the references in \cite{St}), hence
$\overline {\cal B}^3$ does not meet $A(G)$-hull$(K')$. Hence
$A(G)$-hull$(K')$ is contained in a connected component of $\overline
G \setminus \overline {\cal B}^3$ which is impossible.\\

The mentioned results suggested that e.g. Cantor sets in the unit
sphere in $\mathbb C ^2$ are polynomially convex. The purpose of this
paper is to construct Cantor sets in the unit sphere with large
polynomial hull.\\ 

For a subset $A$ of $\mathbb{C}^n$ and a positive number $r$ we
denote by $rA$ the set $\{rz:z\in A\}$. Let $\mathbb{B}^n$ denote the
unit ball in $\mathbb{C}^n$ and $\partial\mathbb{B}^n$ its boundary.\\

\noindent {\bf Theorem.} {\it For any positive number $\beta<1$ there exists a
Cantor set $E$ contained in $\partial\mathbb{B}^2$ whose polynomial hull $\hat{E}$ contains the closed
ball $\overline{\beta\mathbb{B}^2}$. }\\

The theorem implies e.g. that there are continuous arcs in the sphere
(i.e. homeomorphic images of the unit interval of the real line) whose
polynomial hull contains big balls.\\

The basis of the proof is the following observation. Consider two
circles on the sphere, each of them being the intersection of the
sphere with a complex line.  Then the polynomial hull of their union
is equal to the union of their polynomial hulls. However, the
polynomial hull of the union of their $\ep$-neighbourhoods
$(\ep>0)$ is essentially larger than the union of the hull of
their $\ep$-neighbourhoods. More precisely, the following lemma
holds.\\

\noindent {\bf Main lemma.} {\it Let $f$ and $g$ be complex affine functions on
$\mathbb{C}^2$ with $|\nabla f|=|\nabla g|=1$. Suppose the sets
$\{f=0\}\cap\partial\mathbb{B}^2$ and
$\{g=0\}\cap\partial\mathbb{B}^2$ are disjoint circles. Then there
exist positive numbers $a=a(f,g)$ and $r'=r'(f,g)<1$ such that for
any positive $\ep$ the following inclusion holds: \be\label{1}
\{|f\cdot g|\leq a \ep \}\cap (\overline{\mathbb{B}^2}\setminus
r'\mathbb{B}^2) \subset (\{|f|\leq\ep\}\cup
\{|g|\leq\ep\})\cap (\overline{\mathbb{B}^2} \setminus
r'\mathbb{B}^2).  \ee } \\

By a complex affine function $f$ we mean a mapping $f:\mathbb{C}^2\to
\mathbb{C}$ of the form $f(z)=f_0+f_1\cdot z_1+f_2\cdot z_2$ for
complex numbers $f_0, f_1$ and $f_2$. \\

\noindent {\bf Corollary 1.} {\it With $f,\, g,\, a,\, \ep$ and $r'$ as above the
inclusion \be\label{2} \{|fg|\le a\ep\}\cap
\overline{r\mathbb{B}^2} \subset \widehat{\big((\{|f|\leq\ep\}\cap
r\partial\mathbb{B}^2)\cup(\{|g|\leq\ep\} \cap
r\partial\mathbb{B}^2)\big)}\ee holds for any $r\in [r',1]$.}\\

\noindent {\bf Proof.} (1) implies in particular that for any $r\in [r',1]$
\be\label{3} \{|fg|\leq a\ep\} \cap
r\partial\mathbb{B}^2\subset(\{|f|\leq\ep\}\cap
r\partial\mathbb{B}^2)\cup(\{|g|\leq\ep\}\cap r\partial
\mathbb{B}^2) \ee holds. The polynomial hull of the left hand side is
$\{|fg|\leq a\ep\}\cap \overline{r\mathbb{B}^2}$, hence (2)
holds. \hfill $\Box$. \\

\noindent {\bf Remark 1.} {\it Each set on the right hand side of (3)
is the intersection of the sphere $r\partial\mathbb{B}^2$ with the
closed $\ep$-neighbourhood of a circle (i.e. a solid torus in the
3-sphere if $\ep$ is small). Suppose the complex lines $\{f=0\}$ and
$\{g=0\}$ intersect inside $\mathbb{B}^2$. Then the left hand side of
(2) contains e.g. a ball of radius $\sqrt{a}\cdot \sqrt{\ep}$ around
the intersection point. The radius is much larger than $\ep$ for small
enough $\ep>0$. The proof of the theorem is based on this
observation. }\\

\noindent {\bf Proof of the main lemma.} Since the circles are disjoint
there exist $a$ and $r'$ such that also the sets $\{|f|\leq a\}\cap
(\overline{\mathbb{B}^2}\setminus r'\mathbb{B}^2)$ and $\{|g|\leq a\}
\cap(\overline{\mathbb{B}^2}\setminus r'\mathbb{B}^2)$ are
disjoint. If $z\in\overline{\mathbb{B}^2}\setminus r'\mathbb{B}^2$ and
$|f(z)\cdot g(z)|\leq a\cdot\ep$, then either $|f(z)|>a$ and then
$|g(z)|\leq\ep$, or $|f(z)|\leq a$ and then $|g(z)|>a$ and hence
$|f(z)|\leq\ep$. \hfill $\Box$.\\

Note that the constants $a$ and $r'$ can be chosen so that they serve
also for pairs of complex affine functions close to $f$ and $g$.\\

\noindent {\bf Preparation of the proof of the theorem.} We will find the
required set $E$ of the form $E=\bigcap\limits^\infty_{N=1} E_N$ for a
decreasing family of closed sets $E_N\subset\partial\mathbb{B}^2$ such
that $\overline{\beta\mathbb{B}^2}\subset\hat{E}_N$ for each $N$. Then
$\hat{E}\supset \overline{\beta\mathbb{B}^2}$. Indeed, let $z\in
\overline{\beta\mathbb{B}^2}$. For any fixed polynomial $p$ and for
any $N$ there exists a point $z_N\in E_N$ for which $|p(z)|\leq
|p(z_N)|$. If (for fixed $p$ ) $z^*$ is an accumulation point of the $z_N$, then
$z^*\in E$ and by continuity $|p(z)|\leq |p(z^*)|$. This holds for
arbitrary polynomials $p$, hence $z\in\hat{E}$.

Each set $E_N$ will be the finite union of disjoint solid tori of the
form described above. Introduce the following notation. For a complex
affine function $f$ with $|\nabla f|=1$ and a positive number $\sigma$
we denote $T_f(\sigma)=\{z\in\partial\mathbb{B}^2:
|f(z)|\leq\sigma\}$. After a unitary change of coordinates in
$\mathbb{C}^2$ we may assume that $|f|$ has the form $|z_1-(1-s)|$ for
a real number $s=s_f$. The mentioned unitary transformation takes
$T_f(\sigma)$ to the set
$T^s(\sigma)=\{z\in\partial\mathbb{B}^2:|z_1-(1-s)|\leq\sigma\}$. This
is a solid torus if $s<1$ and $\sigma<s$. For $r<1$ denote
$^rT_f(\sigma)=\{|f|\leq \sigma\}\cap
r\partial\mathbb{B}^2$. Similarly we write $^rT^s(\sigma)$. Denote
finally the complex lines (the symmetry axes of the respective tori)
by $\ell_f=\{f=0\}$ and $\ell^s=\{z_1=1-s\}$.

Note that the number $\sigma$ together with the unitary invariant $s$
of a torus  $T_f(\sigma)$ contained in $\partial\mathbb{B}^2$ determine its
diameter. If $s$ and $\sigma$ are small then the diameter of
$T^s(\sigma)$ (hence of each unitarily equivalent torus) is small.\\

By the following easy lemma it is enough to cover bidiscs by
polynomial hulls of suitable Cantor sets. \\

\noindent {\bf Lemma 1.} {\it For any $\beta\in (0,1)$ there is a number
$q\in(0,1)$ such that $\overline{\beta\mathbb{B}^2}$ can be covered by
a finite union of bidiscs of the form $q(\overline{{\cal D}_1}
\times\overline{{\cal D}_2})$ where $\overline{{\cal D}}_j$ are closed
discs in $\mathbb{C}$ centered at zero such that $\partial{\cal
D}_1\times \partial{\cal D}_2\subset\partial\mathbb{B}^2$. \hfill
$\Box$}\\

The following proposition allows to cover bidiscs of lemma 1 by the
polynomial hull of finite unions of disjoint solid tori which are
arbitrarily thin tubular neighbourhoods of (in general not small)
circles.\\

\noindent {\bf Proposition 1.} {\it Let ${\cal D}_1$ and ${\cal D}_2$ be discs
in $\mathbb{C}$ centered at the origin and such that $\partial{\cal
D}_1\times \partial {\cal D}_2 \subset\partial\mathbb{B}^2$. Let
$\gamma>0$ and $q\in(0,1)$. There exist two numbers $s_1,s_2\in(0,1)$
and a positive number $\ep'$, all depending only on $\CD_1,
\CD_2, q$ and $\gamma$, such that for each $\ep\in(0,\ep')$
there exist two families of solid tori
$T_j(\ep)=T_{f_j}(\ep)$ and
$T^*_k(\ep)=T_{g_k}(\ep)$, each family containing finitely
many tori which are unitarily equivalent to $T^{s_1}(\ep)$ and
$T^{s_2}(\ep)$, respectively, with the following properties: 
\begin{itemize}
\item All thicker tori $T_j(2\ep)$ and $T^*_k(2\ep)$ are
pairwise disjoint and contained in the $\gamma$-neighbourhood of
$\partial\CD_1\times\partial\CD_2$.
\item The polynomial hull of the union satisfies the relation 
\end{itemize}

\be\label{4}
q(\overline{\CD_1}\times\overline{\CD_2})\subset
\widehat {\left({\bigcup_j T_j(\ep)\cup\bigcup_k T^*_k(\ep)}\right)}.  \ee }

\noindent {\bf Remark 2.} {\it For some $r' < 1$ depending only on $\CD_1,\CD_2,q$
and $\gamma$ the tori in proposition 1 can be chosen such that for
all $r\in[r',1]$ also the inclusion\\

\noindent $(4_r) \quad\qquad\qquad\qquad\qquad
q(\overline{\CD_1}\times\overline{\CD_2})\subset \widehat{\left({\bigcup_j
\,^rT_j(\ep)\cup \bigcup_k\,^rT^*_k(\ep)}\right)}$\\

\noindent holds. }\\

\noindent Let $\mathbb{D}$ be the unit disc in $\mathbb{C}$ and $\mathbb{T}=\partial\mathbb{D}$.\\

\noindent {\bf Proof of proposition 1.} Let $\CD_1=R_1 \mathbb{D}, \; \CD_2=R_2
\mathbb{D}$. Increasing $q$ we may assume that $q<1$ is as close to
$1$ as needed. Let $\zeta^j_1$ be equidistributed points on $qR_1\TT$
with distance of nearest points being a number between $B\ep$ and
$(B+1)\ep$, where $B$ is any constant, $B\geq 5$, and
$\ep>0$ is small enough. If $\ep$ is small for each $B$ such
points can be found. Similarly, let $\zeta_2^k$ be equidistributed
points on $qR_2\TT$ with distance between closest points in
$[B\ep,(B+1)\ep]$. Define
\[
f_j(z) = z_1-\zeta^j_1, \quad g_k(z)=z_2-\zeta^k_2.
\]
All sets $T_j(2\ep)$ and $T^*_k(2\ep)$ are non-empty tori if
$\ep$ is small enough and they are pairwise disjoint for such
$\ep$. Disjointness is clear for tori of the same family since
$B\ge 5$ and follows from the fact that $\{f_j=0\}\cap\{g_k=0\}$ is
contained in $\BB^2$ in the other case.

For any $j$ and $k$ there is a unitary transformation which takes
$|f_j|$ to $|z_1-(1-s_1)|$,  $ s_1=1-R_1q$, and $|g_k|$ to
$|z_2-(1-s_2)|$, $s_2=1-R_2q$. It follows that the $T_j(\ep)$ are
unitarily equivalent to $T^{s_1}(\ep)$ and the $T^*_k(\ep)$
are unitarily equivalent to $T^{s_2}(\ep)$. Moreover, by
corollary 1 there exist constants $a$ and $r'$ depending only on $s_1$
and $s_2$ (since all pairs $(|f_j|, |g_k|)$ are unitarily equivalent
to $(|z_1-(1-s_1)|, |z_2-(1-s_2)|)\,)$ such that
\[
\{|f_jg_k|\le a\ep\} \cap r\BB^2\subset \widehat{ \big( \left( 
\{|f_j|\leq \ep\} \cup \{|g_k| \leq \ep\} \right) \cap r \partial
\BB^2 \big)}
\]
for all $r\in [r',1]$ and all $j$ and $k$. The left hand side contains
the bidisc $\{|f_j|\leq \sqrt{a\ep}\}\cap
\{|g_k|\leq\sqrt{a\ep}\}=(\zeta^j_1,\zeta_2^k)+\sqrt{a\ep}\;(\DD\times\DD)$
of radius $\sqrt{a}\sqrt{\ep}$ around the intersection point
$\ell_{f_j}\cap\ell_{g_k}=(\zeta^j_1,\zeta^k_2)$ (with $r'$ close to $1$ and
$\ep$ small enough).  If $\sqrt{a}\sqrt{\ep}>(B+1)\ep $ we obtain
(running over all pairs $(j,k)$) that the right hand side of (4$_r$)
contains the product of the circles $qR_1\TT$ and $qR_2\TT$, hence
(4$_r$) holds.

Finally, taking from the beginning $q$ close to $1$ (and $\ep>0$
small), we ensure that all $T_j(2\ep)$ and $T_k^*(2\ep)$ are
contained in a small neighbourhood of
$\partial\CD_1\times\partial\CD_2$. \hfill $\Box$\\

The following proposition allows to cover small enough neighbourhoods
of given complex affine discs (intersections of complex lines with
suitable balls) by polynomial hulls of unions of disjoint tori with
small diameter contained in the sphere. Propositions 2, 3 and 4 below
will be stated for tori related to the function $z_1-(1-s)$ for some
$s\in (0,1)$. They hold for $z_1-(1-s)$ replaced by
any complex affine function $f$ such that $|f|$ is unitarily
equivalent to $|z_1-(1-s)|$.\\

\noindent {\bf Proposition 2.} {\it Let $s\in (0,1)$. For any small $\delta>0$,
any $q\in (0,1)$ close to $1$ and any small $s'>0$ one can find two
positive numbers $s_1<s',\; s_2<s'$ and a positive number
$\ep'=\ep' (s,s_1,s_2,\delta)$ such that for any positive $\ep <\ep'$ there exist two finite
families of complex affine functions $f_j$ and $g_k$ 
with the following properties:\\
\noindent The
related solid tori $T_j(2\ep) \stackrel{{\rm\scriptsize def}}{=}
T_{f_j}(2\ep)$ and $T_k^*(2\ep) \stackrel{{\rm\scriptsize
def}}{=} T_{g_k}(2\ep)$ are all pairwise disjoint and contained
in $\{|z_1-(1-s)|\leq\delta\}$. The $T_j(\ep)$ are unitarily
equivalent to $T^{s_1}(\ep)$ and the $T_k^*(\ep)$ are
unitarily equivalent to $T^{s_2}(\ep)$. Moreover \be\label{5} \overline{q\BB^2}\cap
\{|z_1-(1-s)|\leq\ep\}\subset \widehat{\left( \bigcup_j
T_j(\ep)\cup\bigcup_k T^*_k(\ep)\right) }.  \ee The following stronger
assertion holds. There exists $r'<1$ such that for any $r\in[r',1]$
the inclusion\\

\noindent $(5_r) \quad\qquad\qquad\qquad \overline{q\BB^2}\cap
\{|z_1-(1-s)|\leq\ep\}\subset \widehat{\left( \bigcup_j 
\,^rT_j(\ep)\cup\bigcup_k\,^rT^*_k(\ep)\right) }.$\\

\noindent holds. }\\

Note that the number of tori and their symmetry axes will be chosen
together with $\ep$.\\

\noindent {\bf Remark 3.} {\it Proposition 2 describes a subset of the sphere
whose connected components have a priori given small diameter such
that its polynomial hull contains e.g. the disc $\{z_1=0,\; |z_2|\leq
q\}\quad (q<1)$. This statement has to be opposed to the following
still open
conjecture of Vitushkin which appeared in connection with a problem of
E. Kallin on polynomial convexity of finite unions of disjoint
balls. }\\

\noindent {\bf Conjecture (Vitushkin).} {\it Let $X$ be a relatively closed
complex manifold of dimension one in $\BB^2$ smooth up to the boundary
and transversal to $\partial\BB^2$. Suppose $X$ contains the
origin. Then $X$ has a boundary component of diameter bounded from
below by an absolute constant.}\\

The link between Vitushkin's conjecture and the above statement (or the
theorem, respectively) is the following open problem.\\

\noindent {\bf Problem.} {\it Let $K\subset\mathbb{C}^2$ be a compact,
$K\not=\hat{K}$. Suppose $z \in \hat{K} \setminus K$. Under which conditions for every $\ep>0$ the
$\ep$-neighbourhood $K_{\ep}$ of $K$ contains the boundary of a
Riemann surface passing through $z$?} \\

Note that proposition 2 does not control the ratio of the numbers
$\ep$ (the width of the tubular neighbourhood of $\ell^s$ which is
covered by the polynomial hull of the union of small tori) and
$\delta$ (the width of the tubular neighbourhood which contains the
small tori). Propositions 3 and 4 below will allow such a control, after
which the inductive construction of the $E_N$ can be easily
provided.\\

\noindent {\bf Plan of proof of proposition 2.} The functions $f_j$ in
proposition 2 will be chosen so that $\ell_{f_j}$ will pass through a
point $p_j$ where the $p_j$ are equidistributed on the circle
$\{z_1=1-s\}\cap R\partial\BB^2$ for a suitable $R<1$ and close to $1$,
and the $\ell_{f_j}$ are complex tangent to the sphere
$R\partial\BB^2$. The $\ell_{g_j}$ pass through points $p_j^*$ which
are equidistributed on the same circle. Moreover, the $p_j^*$ are obtained
from  $p_j$ by turning by a fixed angle $\psi$ in $z_2$-direction and the
$\ell_{g_j}$ are obtained by turning the complex tangents to the sphere $R\partial\BB^2$
at $p_j^*$ by a fixed angle $\nu$. Lemma
2 below states that for suitable $p_j, \,R$ and $\ep$ the tubes
$\{|f_j|\leq 2\ep\}\cap\overline{\BB^2}$ are pairwise disjoint. (Hence
the tori $T_j(2\ep)=\{|f_j|\leq 2\ep\} \cap\partial\BB^2$ are pairwise
disjoint and not linked with each other in $\partial\BB^2$).
Moreover, with a suitable choice of the angle $\nu$ the same
is true for the tubes $\{|g_k|\leq 2\ep\}\cap\overline{\BB^2}$ and the
tori $T^*_k(2\ep)$.

On the other hand each complex line $\ell_{f_j}$ intersects several
lines $\ell_{g_k}$ inside $\BB^2$ (equivalently, the corresponding
circles $\ell_{f_j}\cap\partial\BB^2$ and
$\ell_{g_k}\cap\partial\BB^2$ are linked in $\partial\BB^2$). Lemma 3
below allows to choose the angle $\psi$ in such a way that the tubes
$\{|f_j|\le 2\ep\}$ and $\{|g_k|\le 2\ep\}$ intersect at points which
are not contained in the sphere $\partial\BB^2$, i.e. the tori
$T_j(2\ep)$ and $T^*_k(2\ep)$ are disjoint. An application of the main
lemma will give the proposition.\\

\noindent {\bf The tori $T_j(\ep)$ and $T^*_j(\ep)$}. The $T_j(\ep)$ will be
determined by the following parameters: $s\in(0,1)$, a small number
$t>0$, a (large) natural number $\cal N$ and a small number $\ep>0$. The
$T^*_j(\ep)$ will be determined by $s,\, t,\, \cal N,\, \ep$ and a small number
$\nu$.

Consider the intersection
\[
\ell^s\cap\partial\BB^2=\{(1-s,R_2e^{i\phi}); \phi\in [0,2\pi)\} \quad (R_2^2=2s-s^2).
\]
Denote by $C_t$ the slightly smaller concentric circle
$C_t=\{(1-s,R_2(1-t)e^{i\phi}):\phi\in[0,2\pi)\}=\ell^s\cap
R\partial\BB^2$, where $R$ and $t$ are related by the equality
$(2t-t^2)(2s-s^2)=1-R^2$.

Let $p_j$ be equidistributed points on $C_t$,
\[
p_j=(1-s, R_2(1-t)e^{i\phi_j}), \mbox{ where } \phi_j=\frac{2\pi}{\cal
  N} j, \quad j=0,\ldots,{\cal N} -1.
\]
Define the constant $B$ by the following relation: the distance
between nearest of the equidistributed points, $|p_j-p_{j-1}|$, equals
$B\cdot\ep$. Denote by $\ell_j$ the complex lines through $p_j$ which
are tangent to $R\partial\BB^2$, \be\label{6}
\ell_j=\{p_j+v_j\cdot\zeta:\zeta\in\CC\},\quad
v_j=((1-t)R_2,-(1-s)e^{i\phi_j}), \ee and let $f_j$ be annihilating
functions of the $\ell_j$ with gradient of norm one. 

For a (small) number $\psi$ denote by $p_j^*$ the points
\[
p_j^*=\big(1-s, R_2(1-t)e^{i(\phi_j + \psi)}\big).
\]
Finally, for a small constant $\nu>0$ we denote by $\ell^*_j, \;
j=0,\ldots, {\cal N} -1,$ the complex line obtained from the complex
tangent to  $R\partial\BB^2$ at $p_j^*$ by turning by
a fixed angle,  \be\label{7}
\ell^*_j=\{p_j^*+w_j\cdot\zeta:\zeta\in\CC\},\quad
w_j=\left( (1-t)R_2,-(1-s+\nu) e^{i(\phi_j+ \psi)} \right) .  \ee

\noindent Let $g_j$ be annihilating functions of $\ell^*_j$ with gradient of norm one.\\

\noindent {\bf Lemma 2.} {\it a) Let $s\in (0,1)$ and $t>0$ be
sufficiently small. There exist positive constants $\ep'=\ep'(s,t )$
and $B'=B'(s,t)$ such that if $\ep\in(0,\ep')$ and $B>B'$ then for the
$f_j$ defined above for chosen parameters $s,\, t,\, {\cal N}$ and $\ep$
(with $B$ related to ${\cal N}$ and $\ep$ as above) the sets
$\{|f_j|\le 2\ep\}\cap \overline{\BB^2}$ are pairwise disjoint.  \\ b)
 If in addition $\nu$ is small enough (depending on $s$ and $t$) there exist
constants $\ep'=\ep'(s,t,\nu)$ and $B'=B'(s,t,\nu)$ such that if
$\ep\in(0,\ep')$ and $B>B'$ then for the $g_j$ defined for the
parameters $s,\, t,\, {\cal N},\, \nu$ and $\ep$ and an arbitrary
parameter $\psi$ the sets $\{|g_j|\leq
2\ep\}\cap\overline{\BB^2}$ are pairwise disjoint. }\\

\noindent {\bf Lemma 3.} {\it With suitable constants $\ep'$ and $B'$
  which are greater than the constants of lemma 2a and 2b and
  parameters  $s,\, t,\, {\cal N},\, \nu$ and $\ep$ as in lemma 2 with
  $B>B'$ one can choose
the constant $\psi$ in such a way that the tori
$T_j(2\ep)$ and $T^*_k(2\ep)$ of lemma 2 are disjoint for all $j$ and $k$.}\\

\noindent {\bf Proof of lemma 2.} The argument is roughly, that for small
$t>0$ the intersection $\ell_j\cap\overline{\BB^2}$ is a disc of small
diameter. If two such discs intersect their centers $p_j$ must be
close. But then the corresponding complex tangencies to
$R\partial\BB^2$ are ``almost parallel'' so they cannot intersect at
points close to the $p_j$.

More precisely, the set $\{|f_j|\leq 2\ep\}$ is the union of complex
lines ${\cal L}_j=\{\tilde{p}_j+v_j\cdot\zeta:\zeta\in\mathbb C\}$ with the
same direction $v_j$ as $\ell_j$ through points $\tilde{p}_j\in\ell^s$
with distance from $p_j$ not exceeding $A\ep$ for a constant $A$
depending on the angle between $\ell^s$ and $\ell_j$, hence by unitary
equivalence, on $s$ and $t$ only. Hence, \be\label{8} \tilde{p}_j
=(1-s, R_2(1-t)e^{i\phi_j}+\alpha_j) \mbox{ with } |\alpha_j|\leq
A\ep.  \ee Let $j\not= k$. For the intersection point ${\cal L}_j \cap
{\cal L}_k\not= \emptyset$ we have
\[
\tilde{p}_j+v_j\zeta = \tilde{p}_k+v_k\zeta' \mbox{ for some }
\zeta,\zeta'\in\mathbb C.
\]
From (6) and (8) we obtain that $\zeta=\zeta'$ and \be\label{9} \zeta
= \frac{R_2(1-t)}{(1-s)} +
\frac{\alpha_j-\alpha_k}{(1-s)(e^{i\phi_j}-e^{i\phi_k})} \mbox{ with }
|\alpha_j|\leq A\ep, \; |\alpha_k|\leq A\ep.  \ee For small $t>0$ and
$\ep < \ep'(t,s)$ the absolute value $|\zeta|$ of this number can be
estimated from below by a positive constant depending only on $s$,
provided $B\ge B'(s,t)$ and $B'(s,t)$ is chosen so that
\[
|e^{i\phi_j}-e^{i\phi_k}|^{-1} \cdot 2A\ep\leq
 \frac{2A\ep}{B'(s,t)\cdot\ep}< \frac{1}{2} R_2(1-t).
\]
On the other hand $\ell_j\cap\overline{\BB^2}$ is a closed disc of
radius $\sqrt{1-R^2}$ which, for fixed $s$, is small if $t>0$ is
small. Hence for $\ep<\ep' (t,s)$ the diameter of the intersection
${\cal L}_j\cap\overline{\BB^2}$ is small (${\cal L}_j$ as above
contained in $\{|f_j|\leq 2\ep\}$). For those $t$ and $\ep$ the
intersection ${\cal L}_j\cap\overline{\BB^2}$ cannot contain the two
points $\tilde{p}_j$ and $\tilde{p}_j+v_j\zeta={\cal L}_j\cap {\cal L}_k$ since
their distance is bounded from below by a constant depending only on
$s$. Since $\tilde{p}_j\in{\cal L}_j\cap\overline{\BB^2}$ the point
$\tilde{p}_j+v_j\zeta$ is not in $\overline{\BB^2}$. Part a) is
proved.

To prove assertion b) increase $A$ if necessary and replace the number
$(1-s)$ in (\ref {9}) by $(1-s+\nu)$. Use that for $\nu$ small, $\nu<\nu(s,t)$,
and $\ep<\ep'(s,t,\nu)$, the complex lines ${\cal L}^*_j$ parallel to
$\ell^*_j$ and of distance not exceeding $2\ep$ from $\ell^*_j$ still
intersect $\overline{\BB^2}$ along a disc of small diameter. The
remaining arguments are the same as for part a). \hfill $\Box$ \\

\noindent {\bf Remark 4.} {\it Note that the unitary transformation $(z_1,z_2)\to
(z_1,z_2e^{-i\phi_j})$ maps the torus $T_j(\ep)=T_{f_j}(\ep)$ to the
torus $T_0(\ep)=T_{f_0}(\ep)$ and $T^*_j(\ep)=T_{g_j}(\ep)$ to
$T^*_0(\ep)=T_{g_0}(\ep)$. Moreover, the tori $T_j(\ep)$ are unitarily
equivalent to $T^{s_1}(\ep)$ with \be\label{10}
s_1=1-R=(1+R)^{-1}(2t-t^2)(2s-s^2) \ee and $T^*_j(\ep)$ are unitarily
equivalent to $T^{s_2}(\ep)$ for some number $s_2>s_1$ which is close
to $s_1$ if $\nu$ is small.}\\

\noindent {\bf Proof of lemma 3.} We want to choose $\psi$ so that $T_j(2\ep)$ and
$T^*_k(2\ep)$ are disjoint for all $j$ and $k$. Note that the norm
$|P_{j,k}|$ of the intersection point $P_{jk}=\ell_j\cap\ell^*_k$
depends only on $m=k-j$. The idea is the following. When $|P_{0,m}|$
is close to $1$, the
points $P_{0,m}$ form approximately an arithmetic progression with
step  $const \cdot B\ep$ on a real line in the complex line
$\ell_0$. Changing the parameter $\psi$ leads approximately to
translating the approximate arithmetic progression by the parameter
$\psi$ inside the real line. If $B$ is large enough, this allows to choose
$\psi$ in such a way that the intersection points 
$P_{0,m}$ are not in $\partial\BB^2$, moreover, a neighbourhood of
them of size comparable with $\ep$ (containing the intersection
$\{|f_0|\leq 2\ep\}\cap\{|g_m|\leq 2\ep\}$) does not meet $\partial\BB^2$.

Here is the precise argument. Let $s,\,t,\, {\cal N}\, ,\nu$ and
$\ep$ be chosen according to
lemma 2 with constants $\ep'$ and $B'$ specified
below and greater than the constants in part a and b of the lemma.
We will change the parameter $\psi$. The points $p_k^*$, the complex lines
$\ell^*_k$, the intersection points  $P_{jk}=\ell_j\cap\ell^*_k$ and the tori $T^*_k(\ep)$ will depend on $\psi$. We will indicate the
dependence on $\psi$ only sometimes when we want to draw special
attention to this fact.

From (\ref{6}) and (\ref{7}) it follows that the intersection point
$P_{j,k}=P_{j,k}(\psi)$ of the complex lines $\ell_j$
and $\ell^*_k$ is determined by
\[
P_{j,k}(\psi)=p_j+v_j\cdot\zeta(\psi)=p_k^*(\psi)+w_k(\psi)\cdot\zeta'(\psi)
\]
for some $\zeta(\psi),\zeta'(\psi)\in {\mathbb C}$. From the same formulas
we obtain that $\zeta(\psi)=\zeta'(\psi)  \stackrel{{\rm\scriptsize def}}{=}  \zeta_{j,k}(\psi)$ and \be\label{11}
\zeta_{j,k}(\psi)=\frac{R_2\cdot(1-t)}{1-s}\,
\frac{1-e^{i(\phi_k-\phi_j+\psi)}}{1-Q e^{i(\phi_k-\phi_j+\psi)}}, \ee where
\[
Q=(1-s+\nu)(1-s)^{-1}>1.
\]
Here as above we put $\phi_j=\frac{2\pi j}{\cal N}$ and assume that
the natural number $\cal N$ is big. Put \be\label{12} F(\phi)=\Big|
\frac{1-e^{i\phi}}{1-Q e^{i\phi}}\Big| = \Big|
\frac{2\sin\frac{\phi}{2}}{1-Q e^{i\phi}}\Big|, \ee so that
\be\label{13} |\zeta_{j,k}(\psi)| = \frac{R_2(1-t)}{1-s}
F(\phi_k-\phi_j+\psi).  \ee Note that the function $(|1-Q e^{i\phi}|)^{-1}$
is of class $C^\infty(\mathbb{R})$, hence for small $\phi\not= 0$ we have
\[
F'(\phi)=\frac{\cos \frac{\phi}{2}\cdot \mbox{sgn}\phi}{|1-Q
e^{i\phi}|} + |2\sin\frac{\phi}{2} | \cdot \Big(\frac{1}{|1-Q
e^{i\phi}|}\Big)'.
\]
Hence \be\label{14} 0< C_Q \le |F'(\phi)|\le 2 C_Q \mbox{ for } \phi\not=0 \mbox{
  and } |\phi| \leq \Phi_Q \ee for  positive
  constants  $\Phi_Q$ and $C_Q$ which may be chosen depending only on $Q$. Since $p_j\in R\partial\BB^2$ and
  $v_j$ is the direction of the complex tangent line to
  $R\partial\BB^2$ we have
\[
|P_{j,k}|^2 = |p_j+v_j\cdot\zeta_{j,k}|^2 = |p_j|^2 + |v_j|^2 |\zeta_{j,k}|^2.
\]
By (\ref{6})    $\;\;|v_j|^2=|p_j|^2=R^2$. Hence 
\be\label{15}
|P_{j,k}(\psi)|^2 = R^2 \left(1+ \frac{R^2_2(1-t)^2}{(1-s)^2} F^2(\phi_k-\phi_j+\psi)\right),
\ee
with 
\be\label{16}
R^2=1-(2t-t^2)(2s-s^2).
\ee

The set $\{|f_j|\leq 2\ep\}\cap\{|g_k|\leq 2\ep\}$ is
contained in the (closed) $A'\ep$-neighbourhood of $P_{j,k}=P_{j,k}(\psi)$. The
constant $A'$ depends only on $\nu,\;s$ and $t$. For the square of the
norm of points in $\{|f_j|\le 2\ep\} \cap \{|g_k|\leq 2\ep\}$ this
implies that these numbers are contained in the (open) $4A'\ep$-neighbourhood
of $|P_{j,k}(\psi)|^2$ provided $|P_{j,k}(\psi)|\leq \frac{3}{2}$ and
$A'\ep < 1$.

We want to choose $\psi$ so that the $4A'\ep$-neighbourhoods of the
$|P_{j,k}(\psi)|^2$ do not contain $1$. Put $m=k-j$. Then $\phi_k-\phi_j+\psi =
\phi_m+\psi =- \phi_{-m}+\psi  $ and $|P_{j,k}(\psi)|^2 = |P_{0,m}(\psi)|^2 = |P_{0,-m}(-\psi)|^2$. Hence all possible
values of (\ref{15}) are obtained when $j=0$ and $k=m$ runs over
integers.

If for some $\psi$ and $m_0$ 
 \be\label{17}
\Big| |P_{0,m_0}(\psi)|^2 -1\Big| < 4A'\ep,  \ee 
then by (\ref{15}) and (\ref{16}) $\;\;F^2(\phi_{m_0}+\psi)\; t^{-1}$ can
be estimated from above and from below by positive constants depending
only on $s$ provided $\ep$ is small compared with $t$. 
Hence in this case  $\phi_{m_0}+\psi$ is comparable to either
$\sqrt{t}$ or $-\sqrt{t}$.

Suppose now that for for some $m_0$ the $4A'\ep$-neighbourhood of
$|P_{0,m_0}(0)|^2  $ contains $1$. (Otherwise we are done.) Note that
by symmetry this is true also for $-m_0$. If $t$ is small then by the
preceding arguments  $\phi_{m_0}$ and  $\phi_{-m_0}$ are much smaller
than $\Phi_Q$. By (\ref{14}) and (\ref{15}) for fixed $\psi$
close to $0$ and  for integer numbers
$m$ close to $m_0$ (close to $-m_0$, respectively) the possible values of the right hand side of
(\ref{15}) have distance from each other at least equal to $c(s,t)\cdot R_2(1-t)|\phi_{m+1}-\phi_m|$ for a positive
constant $c(s,t)$ depending on $s$ and $t$. Since
$|\phi_{m+1}-\phi_m|=\frac{2\pi}{\cal N}>2 \sin\frac{\pi}{\cal N}$ and
$R_2(1-t)2\sin\frac{\pi}{\cal N}=B\ep$, the difference of the values
of (\ref{15}) for fixed $\psi$ close to $0$ and  for $m$ close to
$m_0$ (close to $-m_0$, respectively)  is at
least $c(s,t)\cdot B \cdot\ep$. Take the constant
$B'$ large enough so that the
latter constant is at least  $40A'\ep$.

Using again (\ref{14}) and (\ref{15}) take $\psi$ so that 

\[
10A'\ep < ||P_{0,{m_0}}(\psi)|^2 - |P_{0,{m_0}}(0)|^2 | < 30A'\ep
\]

{\noindent}and the same estimate holds for  $m_0$ replaced by $-m_0$. Note that
$\psi$ is comparable to $\ep$ with multiplicative constants depending
on $s$ and $t$ and, hence, $\psi$ is  small if $\ep$ is small. The $4A'\ep$-neighbourhood of
$|P_{0,m_0}(\psi)|^2  $ does not contains $1$ and the same is true for $m_0$
replaced by $-m_0$. The above arguments give that for all $m$ close to $m_{0}$
\[
||P_{0,m}(\psi)|^2-1|> 4A'\ep,
\]

\noindent hence, since the function $F$ is strictly monotonic on
the positive half-axis (respectively, on the negative half-axis) this
holds for all $m$.
 
With this choice of $\psi$ we obtained that for all $j$ and $k$ the set $\{|f_j|\leq
2\ep\}\cap\{|g_k|\leq 2\ep\}$ does not meet $\partial\BB^2$ 
if $\ep$
is small enough, hence the tori $T_j(2\ep)$ and $T^*_k(2\ep)$  are disjoint. \hfill
$\Box$\\

\noindent {\bf Proof of proposition 2.} For $s\in(0,1)$ as in the statement we
first choose $t$ and $\nu$ small enough, so that the numbers $s_1$ and
$s_2$ are less than $s'$ (see (\ref{10}) and remark 4) and $R(t)>q$
(see (\ref{16})).  We may also achieve by choosing $t$ and $\nu$ small
that with any choice of the natural number $\cal N$, of the parameter
$\psi$ and of the small
enough positive number $\ep$ the sets
$\{|f_j|\leq 2\ep\}$ and $\{|g_k|\leq 2\ep\}$ are contained in
$\{|z_1-(1-s)|\leq\delta\}$. Choose by lemmas 2 and 3 the relation of the numbers 
$\cal N$, $\ep$ and $\psi$ so that
 the tori $T_j(2\ep)$ and
$T^*_k(2\ep)$ are pairwise
disjoint. They are unitarily equivalent to $T^{s_1}(\ep)$ and
$T^{s_2}(\ep)$, respectively.

It remains to prove (\ref{5}). Apply the main lemma to the pair
$\ell_j$ and $\ell_j^* = \ell_j^*(\psi)$ for any $j$. Since the pairs obtained for
different $j$ are unitarily equivalent to $\ell_0$ and
$\ell_0^*(\psi)$ (see remark 4) and $\ell_0^*(\psi)$ is close to $\ell_0^*(0)$
if $\psi$ is small, there exist
numbers $a>0$ and $r'\in (0,1)$ depending only on $s,\, t$ and
$\nu$  such that (\ref{1})
holds with $f$ replaced by $f_j$ and $g$ replaced by $g_j$. By
corollary 1 for $r\in [r',1]$ the polynomial hull of $\;^rT_j(\ep)\cup
\;^rT_j^*(\ep)$ contains
\[
\{|f_jg_j| \leq a\ep\}
\cap\overline{r\BB^2}\supset\{|f_j|\leq\sqrt{a\ep}\}\cap\{|g_j|\leq
\sqrt{a\ep}\}\cap \overline{r\BB^2}.
\]

Since $f_j(p_j)=0$ and $g_j(p_j^*)=0$, we obtain
$|g_j(p_j)|=|g_j(p_j)-g_j(p_j^*)|\le |p_j-p_j^*| = O(\ep)$. This
implies that for small $\ep$ the latter set contains the bidisc
$p_j+b\sqrt{a\ep}\;(\overline{\DD}\times\overline{\DD})$ of radius
$b\cdot\sqrt{a\ep}$ around the point $p_j$ for some constant $b$
depending only on $s$, $t$ and $\nu$. Recall that the $p_j$ are
equidistributed on the circle $C_t=\{1-s\}\times \{|z|=R_2(1-t)\}$
with distance between nearest points equal to $B\ep$. If
$b\cdot\sqrt{a\ep}>B\ep$, i.e. $\ep<\frac{b^2}{B^2}a$, then the
polynomial hull of $\bigcup\limits_j (T_j(\ep)\cup T^*_j(\ep))$
contains $\{z_1\}\times \{ |z| = R_2(1-t)\}$ for all $z_1$ with
$|z_1-(1-s)|\leq B\ep$. Since $(1-s)^2+(R_2(1-t))^2=R(t)^2>q^2$ we
obtain that for
$\ep<\ep'$ with a suitable choice of the constants $r'$ and
$\ep'(s,s_1,s_2,\delta)$  the inclusion
(5$_r$) holds for $r\in[r',1]$. The weaker inclusion (\ref{5}) follows. \hfill $\Box$\\

The proof of proposition 2 does not give good estimates for the ratio
of the constants $\ep$ and $\delta$.  The point is that the main lemma
gives useful effects essentially only for the pair $\ell_j$ and
$\ell^*_j$ of complex lines but not for arbitrary pairs $\ell_j$ and
$\ell^*_k$ when the intersection point may be close to $\partial\BB^2$
and the constant $a$ of the main lemma is not bounded away from
zero. Below we will state and prove the stronger proposition 4 which
can be directly used for inductive construction of the $E_N$.

The first step towards this goal is proposition 3 below which is in
the spirit of proposition 1. For its proof we use again two families 
of complex lines. The first family,
$\ell^*_k$, consists of certain complex lines parallel to the
$z_2$-axis. For some small $\sigma>0$ and a suitable constant $\alpha$
their $z_1$-coordinates form an $\alpha$-net of the
$\frac{4}{3}\sigma$-neighbourhood of $1-s$ ($s\in(0,1)$). Note that
the intersection of these lines with $\partial\BB^2$ are circles of
diameters comparable with $\sqrt{s}$. The second family, $\ell_j$,
consists of complex tangent lines to a smaller sphere through
equidistributed points on a circle $C_t$ contained in
$\ell^s=\{z\in{\mathbb C}^2 :z_1=1-s\}$. This time $t$ cannot be
chosen arbitrarily small. It will take a value which is determined by
$s$ and $\sigma$. As in proposition 1 the main lemma will be applied
to all pairs $(\ell_j,\ell_k^*)$.\\

\noindent {\bf Proposition 3.} {\it Let $s\in(0,1)$. There exist
positive constants  $\sigma'=\sigma'(s)$
and $\alpha=\alpha (s,\sigma')$  such that the
following holds.\\ Fix
$\sigma\in(0,\sigma')$. Let $\zeta_k$ be all points of the disc
$\{\zeta\in{\mathbb C}:|\zeta-s|\leq\frac{4}{3}\sigma\}$ which are
contained in the lattice $s + \alpha
\mathbb{Z}+i\alpha\mathbb{Z}$. Denote $g_k(z)=z_1-(1-\zeta_k)$.\\
There
exist positive constants $\ep'=\ep'(s,\sigma)$ and
$s^*=s^*(s,\sigma)<s,\; s^*(s,\sigma)\to 0$ for $\sigma\to 0$,
such that for any $\ep\in
(0,\ep')$ the following is true.\\
One can find finitely many complex affine functions $f_j$
such that the tori $T_j(\ep)=\{z\in\partial\BB^2:|f_j(z)|\leq\ep\}$
are unitarily equivalent to $T^{s^*}(\ep)$. Moreover, the tori
$T_j(2\ep)$ and $T^*_k(2\ep)=\{z\in\partial\BB^2;|g_k(z)|\leq 2\ep\}$
are all pairwise disjoint and contained in $\{|z_1-(1-s)|\leq
2\sigma\}$. Finally, for the polynomial hull of the union of the tori
we have \be\label{18} \{|z_1-(1-s)|\leq\sigma\}\cap
(1-\rho)\overline{\BB^2}\subset \widehat {\left(\bigcup_j
T_j(\ep)\cup\bigcup_k T^*_k(\ep)\right)} \ee for a constant
$\rho=\rho(s,\sigma)$ such that $\rho(s,\sigma)\to 0$ for $s$ fixed
and $\sigma\to 0$.\\ The following sharper statement holds.\\ There
exists a constant $r'=r'(s,\sigma)>0$ such that for $r\in[r',1]$\\

\noindent (18$_r$) $\qquad\qquad
\{|z_1-(1-s)|\leq\sigma\}\cap(1-\rho)\overline{\BB^2}\subset
 \widehat{\left(\bigcup_j\,^rT_j(\ep)\cup\bigcup_k\,^rT^*_k(\ep)\right)}.$ }\\

\noindent {\bf Proof.} Let $s$ and $\sigma$ be given as in the statement of the
proposition. The complex lines $\ell_j$ will be obtained as in the
plan of proof of proposition 2 for some parameters $t$ and $\cal N$:
they will be complex tangent to the sphere $R(t)\partial\BB^2$ (see
(\ref{16})) through equidistributed points $p_j\in C_t$ with distance
$|p_{j+1}-p_j|$ depending on $\cal N$.

Determine now the parameter $t$ (so far $\cal N$ is arbitrary). Let
$t$ be maximal so that
$\ell_j\cap\partial\BB^2\subset\{|z_1-(1-s)|\leq\frac{5}{3}\sigma\}$
and, hence, for this parameter $t\;\;$
$\ell_j\cap\partial\BB^2\subset\{|z_1-(1-s)| = \frac{5}{3}\sigma\}$
(see (\ref{6})).
This number $t$ does not depend on $\cal N$ and
$j$ and tends to zero for $\sigma\to 0$. Hence, also $s^*(s,\sigma)\to 0$ for
$\sigma\to 0$ . Moreover, if $\sigma'$
is small, then there is a uniform estimate from below of the angle of
intersection of $\ell^{s'}$ and $\ell_j$ for $s'\in(s-2\sigma,s+2\sigma),\;\;
\sigma\in(0,\sigma')$, $t$ related to $\sigma$ as just described. Indeed, $\ell^{s'}$ is
transversal to the complex tangent space of $\partial\BB^2$ with a
uniform estimate of the angle for the mentioned $s'$. Hence,
if $\sigma'$ is small and $\sigma<\sigma'$ then $t$ is small, hence $R(t)$
is close to one, and,
hence, the $\ell_j$ (being complex tangents to $R(t)\partial
\mathbb{B}^2$ ) are transversal to the $\ell^{s'}$ with uniform estimate
of the angle.

The positive number $\alpha$ will be specified later. Let
$\zeta_k,\; g_k$ and $\ell^*_k = \{g_k = 0\}$ be as in the statement of the
proposition. The $\ell^*_k$ are parallel to $\ell^s$ and
$\ell^*_k\cap\partial\BB^2$ is contained in
$\{|z_1-(1-s)|\leq\frac{4}{3}\sigma\}$ while $\ell_j\cap\partial\BB^2$
is contained in $\{|z_1-(1-s)|=\frac{5}{3}\sigma\}$. The above
observations imply two facts.

First, if $\sigma'$ and hence $t$ are small
and $\ep$ is small depending on $s'$ and $\sigma$ and, moreover, the
parameters $\cal N$ and $\ep$ satisfy the conditions of lemma 2a, then
all tori $T_j(2\ep)$ and $T^*_k(2\ep)$ are disjoint and contained in
$\{|z_1-(1-s)|\le 2\sigma\}$.

Secondly, the main lemma can be applied to all pairs $(\ell_j,\ell^*_k)$ with
uniform constants $a$ and $r'$.

Moreover, using the notation $ \ell^\xi=\{z\in \mathbb C ^2:z_1=1-\xi\}$ with
$\xi\in\mathbb C$, $|\xi-s|\leq\frac{4}{3}\sigma$, we obtain that for each
positive number $\eta$ the set $\ell^\xi\cap\{|f_j|\leq\eta\}$ is a closed
disc in $\ell^\xi$ of radius at least $D\cdot \eta$ around the intersection
point $\ell^\xi\cap\ell_j$ with constant $D$ depending only on $s$ and
$\sigma'$.

Prove now the assertion on the polynomial hull of the union of the
tori for suitably chosen $\alpha$. Note first that for
$|\xi-s|\leq\frac{4}{3}\sigma$ the intersection points
$\ell^\xi\cap\ell_j$,  $j=0,\ldots,{\cal N}-1$, are equidistributed on the
circle $\{1-\xi\}\times\{|z|=R'_\xi\}\subset\BB^2$ for a number
$R'_\xi$ which is close to $\sqrt{1-|1-\xi|^2}$ if $\sigma<\sigma'$
is small. This follows from the fact that the unitary transformation
$(z_1,z_2)\to (z_1,e^{i\phi_j}z_2)$ maps the pair $(\ell^\xi,\ell_0)$
to the pair $(\ell^\xi,\ell_j)$.

Corollary 1 implies that for $r\in[r',1]$ and all $k$ and $j$ the
polynomial hull of $^rT^*_k(\ep)\cup\; ^rT_j(\ep)$ contains the set
\be\label{19} \{|z_1-(1-\zeta_k)|\cdot|f_j(z)|\leq a\ep\} \cap
r\overline{\BB^2}.  \ee Let $\ep$ be as small as required above,
i.e. $\ep<\ep'(s,\sigma)$, and let $B>B'$ with $B'$ the constant of
lemma 2a. Choose $\cal N$ so that the distance between nearest points
$|p_{j+1}-p_j|$ is between $B\ep$ and $(B+1)\ep$. If $\sigma$ is small
then for $|\xi-s|\leq\frac{4}{3}\sigma$ the distance between nearest
of intersection points $\ell^\xi\cap\ell_j$,  $j=0,\ldots,{\cal N}-1$,
does not exceed $2(B+1)\ep$. Take a constant $C$
so that $C\cdot D>2(B+1)$. The set in (\ref{19}) contains the set
\[
\{z\in\overline{r\mathbb B^2}:|z_1-(1-\zeta_k)|\leq\frac{a}{C},\; |f_j(z)|\leq C\ep\}.
\]
Hence, for any fixed $\xi$ with $|\xi-\zeta_k|\leq \frac{a}{C}$ and
for $r$ as above the polynomial hull of $^rT^*_k(\ep)\cup \;
^r\bigcup_j T_j(\ep)$ contains the circle $\{1-\xi\}\times
\{|z|=R'_\xi\}$. Take any $\alpha <\frac{a}{C}$. Then for
$r\in[r',1]$
\[
\widehat{\left(\bigcup_k\; ^rT^*_k(\ep)\cup\bigcup_j\; ^rT_j(\ep)\right)}\supset
\bigcup_{|\xi-s|\leq\sigma} \{1-\xi\}\times \{|z| \leq R'_\xi\}.
\]
The right hand side contains
$\{|z_1-(1-s)|\leq\sigma\}\cap(1-\rho)\overline{\BB^2}$ for suitable
$\rho=\rho(s,\sigma)$ with $\rho(s,\sigma)\to 0$ for $\sigma\to
0$. \hfill $\Box$\\

Now we are ready to state and prove the main proposition.\\

\noindent {\bf Proposition 4.} {\it Let $s\in (0,1)$. There exists a positive
constant $\sigma'=\sigma'(s)$ such that for any $\sigma\in(0,\sigma')$ there exist
finitely many numbers $s_m\in(0,1)$ (their number depends on $s$ and
$\sigma$ and each of them tends to zero for $\sigma\to 0$ ) and a
positive number $\ep'=\ep'(s,\sigma,s_m)$ such that the following
holds. \\ For any $\ep\in(0,\ep')$ there exist
finitely many complex affine functions $f_n$ and related tori
$T_n(2\ep)$, each of them unitarily equivalent to $T^{s_m}(2\ep)$ for
some $m$ with $T_n(2\ep)$ pairwise disjoint and contained in
$\{|z_1-(1-s)|\leq 2\sigma\}$ and such that \be\label{20}
(1-\rho)\overline{\BB^2}\cap\{|z_1-(1-s)|\leq \sigma\}\subset
\widehat{\left(\bigcup_n T_n(\ep)\right)} \ee for a constant $\rho=\rho(s,\sigma)$
such that for fixed $s$ $\rho(s,\sigma)\to 0$ if $\sigma\to 0$. Moreover,
there exists $r'=r'(s,\sigma,s_m)$ such that for $r\in [r',1]$\\ 

\noindent {\rm (20$_r$)} \qquad\qquad\qquad
$(1-\rho)\overline{\BB^2}\cap\{|z_1-(1-s)|\leq\sigma\}\subset
\widehat{\left(\bigcup_n\, ^rT_n(\ep)\right)}$.\\ }

\noindent {\bf Proof.} Let $\sigma',\, r'$ and
$\alpha$ be as in proposition 3, let $\sigma <
\sigma'$ and let  $\ep'(s,\sigma)$ be the constant from the statement
of proposition 3. Proposition 3 gives two families of functions $f_j$
and $g_k$. If $\sigma$ is small the tori related to the $f_j$ have
small diameter, but the tori related to the $g_k$ have large diameter.
Our aim is to apply proposition 2 to each $g_k$.

Let $s(g_k)$ be the number for which $|g_k|$ is unitarily equivalent to
$|z_1-(1-s(g_k))|$. Apply for each $k$ proposition 2 with $s=s(g_k),\,
\delta=\frac{\alpha}{3}$ and $q=r'$. For each $k$ we obtain two
numbers $s_{k,1}, s_{k,2}\in (0,1)$ (it can be achieved that they are
as small as we wish) and a bound for $\ep$ below of which we can find
finitely many non-intersecting tori $T(2\ep)$, unitarily equivalent to
$T^{s_{k,1}}(2\ep)$ or $T^{s_{k,2}}(2\ep)$, contained in
$\{|z_1-\left( 1-s(g_k)\right) |\leq \frac{\alpha}{3}\}$ with the following
property\\ 
\setcounter{alteqn}{\value{equation}}
\setcounter{equation}{0} \addtocounter{alteqn}{1}

\noindent {\rm (21$_r$)} \qquad\qquad\qquad 
 $ \overline{r'\BB^2}
\cap\{|z_1-(1-s(g_k))|\leq\ep\}\subset \widehat{\left(\bigcup\;
^rT(\ep)\right)} $\\

\setcounter{equation}{\value{alteqn}}
\renewcommand{\theequation}{\mbox{\mbox{\arabic{equation}}}} 

\noindent for all $r\in[\tilde{r}_k,1]$ with a suitable $\tilde{r}_k$.

Assume that  $\ep$ is less than all the mentioned
bounds and also $\ep < \ep'(s,\sigma)$ . Consider the second family of tori  $T_j(2\ep)$ from
proposition 3. The $T_j(2\ep)$ are related to complex affine functions 
$f_j$ and to complex lines  $\ell_j = \{f_j=0\}$. 

Labeling the collection of all tori above, i.e. the tori obtained for
each $k$ by proposition 2 and the tori $T_j(2\ep)$ from proposition 3, we obtain the
family $T_n(2\ep)$. The $T_n(2\ep)$ are pairwise disjoint. Indeed by the
choice of $\delta$ the tori obtained by proposition 2 for different
$k$ do not intersect. Since
$\ell_j\cap\partial\BB^2\subset\{|z_1-(1-s)|=\frac{5}{3}\sigma\}$ and
$\{g_k=0\}\cap\partial\BB^2\subset\{|z_1-(1-s)|\leq
\frac{4}{3}\sigma\}$ the tori obtained by proposition 2 do not meet
the $T_j(2\ep)$ from proposition 3, if $\ep$ is small. The inclusion (18$_{r'}$) from
proposition 3 together with (21$_r$), the latter applied for each $k$ with $r\in [\tilde{r}_k,1]$,  imply (20$_r$) with $r$ larger than the maximum of all
$\tilde{r_k}$ and also larger than $r'$ from proposition 3. Proposition 4 is
proved. \hfill $\Box$\\

Proposition 4 enables us to construct the sets $E_N$ inductively.\\

\noindent {\bf Proof of the theorem. \\ Step 1.} Lemma 1 and
proposition 1 give us for any sufficiently small $\ep>0$ a finite
collection of tori. More precisely, we obtain a number $r_1$ (see
remark 2) and a finite collection of numbers $s^{(1)}_m$ such that for
each small $\ep>0$ for each $m$ there exists a finite number of tori
unitarily equivalent to $T^{s^{(1)}_m}(\ep)$ with the following
properties. Denote all the tori by $T^{(1)}$ skipping labeling indices
but indicating that they are obtained at step 1. The $T^{(1)}(2\ep)$
are pairwise disjoint and \be\label{22} \overline {\beta{\BB^2}}
\subset\widehat{\left(\bigcup \, ^{r_1}T^{(1)}(\ep) \right)}.  \ee The
number $\ep$ will be chosen at the second step of the induction. (Note
that the number of tori and the choice of the complex lines which are
the symmetry axes of the tori also depend on $\ep$.)\\

\noindent {\bf Step 2.} For any $s^{(1)}_m$ of step 1 we apply
proposition 4. We obtain a bound
$\sigma'(s^{(1)}_m)$. Choose now the number $\ep$ of step 1 (and hence
the tori of that step). It has to satisfy the following
requirements. First it has to be so small that it fits for step
1. Next for each $m$ we require $\ep<\sigma'(s^{(1)}_m)$. Further, for
chosen $s=s^{(1)}_m$ and $\sigma=\ep$ proposition 4 allows to find
new $s$-parameters which tend to zero for $\sigma=\ep\to 0$. We denote
the collection (over all $m$) of all new $s$-parameters at the second
step by $s^{(2)}$ (we skip indices labeling them) and require
$\sigma=\ep$ being so small that each $s^{(2)}$ is less than
$\frac{1}{2^2}$. Finally proposition 4 asserts for each $m$ and $\ep$
the existence of a number $\rho(s^{(1)}_m,\ep)$ (see (\ref{20})). We
wish the estimate $1-\rho(s^{(1)}_m,\ep)>r_1$ for each $m$. Choose an
$\ep>0$ satisfying all these conditions and denote it by $\ep_1$.

Denote $E_1=\bigcup T^{(1)}(2\ep_1)$. $E_1$ is the union of pairwise
disjoint closed tori of diameter determined by the $s^{(1)}$ and
$\ep_1$. (\ref{22}) implies that
\[
\overline{\beta\BB^2}\subset\hat{E}_1.
\]
Indeed, $T^{(1)}(2\ep_1)\supset T^{(1)}(\ep_1)$ and
$\widehat{T^{(1)}(\ep_1)}\supset\; ^{r_1}T^{(1)}(\ep_1)$ for any of
the tori $T^{(1)}(\ep_1)$.

The set $E_1$ is constructed. Describe now the construction of the set
$E_2$ modulo the choice of the parameter $\ep_2$.

Denote the functions corresponding to the $T^{(1)}(\ep_1)$ by
$f^{(1)}$ (omitting indices as above). From applying proposition 4 to
each of the $T^{(1)}(\ep_1)$ we got the collection of mentioned above
numbers $s^{(2)}<\frac{1}{2^2}$ and we obtain a number $r_2$ (see (20$_r$)) such that the following holds: For
each sufficiently small $\ep>0$ there exist
finitely many pairwise disjoint tori $T^{(2)}(2\ep)$ contained in
$E_1$, each of them unitarily equivalent to a torus
$T^{s^{(2)}_m}(2\ep)$ for some of the numbers $s^{(2)}_m$ of the
collection $s^{(2)}$, such that (since $1-\rho(s^{(1)}_m,\ep)>r_1$ for
each $m$)
\[
\overline{r_1\BB^2}\cap \bigcup\{| f^{(1)}|\leq\ep_1\}\subset
\widehat{\left(\bigcup\; ^{r_2}T^{(2)}(\ep)\right)},
\]
hence by (\ref{22})
\[
\overline{\beta{\BB^2}}\subset \widehat{\left(\bigcup\; ^{r_2}T^{(2)}(\ep)\right)}.
\]
This describes the construction of the set $E_2$ modulo the choice of
the parameter $\ep_2$.\\

\noindent {\bf Step N.} This is the general step of the induction. Let
$N>2$. Suppose the sets $E_1\supset \ldots \supset E_{N-2}$ are
constructed with $\hat{E}_{N-2}\supset\overline{\beta{\BB^2}}$ and the
construction of the set $E_{N-1}$ is described modulo the choice of
the parameter $\ep_{N-1}$.  We want to choose the parameter
$\ep_{N-1}$ and describe the construction of $E_N$ modulo the choice
of the $\ep_N$. 

More precisely, we suppose the following has been
done at step $N-1$. There was found a finite collection of numbers
$s^{(N-1)}$, all less than $\frac{1}{2^{N-1}}$, and a number
$r_{N-1}\in (0,1)$ such that for each sufficiently small $\ep>0$ there
exist finitely many pairwise disjoint tori $T^{(N-1)}(2\ep)$
contained in $E_{N-2}$, each of them unitarily equivalent to a torus
$T^{s^{(N-1)}_m}(2\ep)$ for some number $s^{(N-1)}_m$ of the
collection $s^{(N-1)}$, with the following relation for the polynomial
hull \be\label{23} \overline{\beta\BB^2} \subset
\widehat{\left(\bigcup\; ^{r_{N-1}}T^{(N-1)}(\ep)\right)}.  \ee

We want to choose a suitable number $\ep_{N-1}$ for $\ep$.  Here are the
requirements for the number $\ep$. First it has to be small enough for
step $N-1$ to go through. Further, apply proposition 4 with $s$ to be any of
the numbers $s^{(N-1)}$. We obtain an upper bound for
$\sigma$. Require that $\ep<\sigma'$ the least of all such upper
bounds over numbers $s^{(N-1)}$. Now for $s$ equal to any of the
$s^{(N-1)}$ and $\sigma=\ep$ proposition 4 allows to choose new
$s$-parameters which tend to zero together with $\sigma=\ep$. Require
$\ep$ being so small that all the new $s$-parameters denoted by
$s^{(N)}$ are less than $2^{-N}$. Finally, the proposition asserts the
existence of a number $\rho(s^{(N-1)}_m,\ep)$ for each number
$s^{(N-1)}_m$ of the $s^{(N-1)}$ (see (\ref{20})). We wish the
estimate $1-\rho(s^{(N-1)}_m,\ep)>r_{N-1}$ for each $m$. Choose a
number $\ep>0$ satisfying all these conditions and denote it by
$\ep_{N-1}$. Put
\[
E_{N-1}=\bigcup T^{(N-1)}(2\ep_{N-1}).
\]

Thus the set $E_{N-1}$ is constructed (recall that also the number of
the tori and the choice of their symmetry axes depend on $\ep_{N-1}$)
and by (\ref{23})
\[
\overline{\beta\BB^2}\subset\hat{E}_{N-1}.
\]

Use now proposition 4 to define $E_N$ modulo the choice of $\ep_N$. Denote the
functions corresponding to the tori $T^{(N-1)}(\ep_{N-1})$ of the previous
generation by $f^{(N-1)}$.  We already mentioned the numbers $s^{(N)}$
obtained by proposition 4. Further this proposition (applied with $s$ being
any of the $s^{(N-1)}$ and with $\sigma = \ep_{N-1}$) gives for each
sufficiently small $\ep>0$  finitely many
pairwise disjoint tori $T^{(N)}(2\ep)$ contained in $E_{N-1}$ (each of them
unitarily equivalent to $T^{s^{N}_m}(2\ep)$ with $s^{N}_m$ being
some of the numbers $s^{N}$) and a number $r_N \in (0,1)$ with the
following relation for the polynomial hull \be\label{24}
\overline{r_{N-1}\BB^2}\cap \bigcup\{| f^{(N-1)}|\leq\ep_{N-1}\}\subset
\widehat{\left(\bigcup\; ^{r_{N}}T^{(N)}(\ep)\right)}. \ee

Here we used that $1-\rho(s^{(N-1)}_m,\ep)>r_{N-1}$ for each $m$.
Hence by (\ref{23}) with the $\ep$ in (\ref{23}) replaced by
$\ep_{N-1}$ we have
\[
\overline{\beta\BB^2} \subset \widehat{\left(\bigcup\;
^{r_{N}}T^{(N)}(\ep)\right)}
\]
for sufficiently small positive $\ep$ of step $N$ and the tori $T^{(N)}(\ep)$ the
existence of which is obtained at step $N$.\\

The induction is complete, hence the theorem is proved. \hfill $\Box$\\

After the paper was written L.Stout informed us that he constructed a
Cantor set in $\mathbb{C}^n$ with non-trivial polynomial hull in case
$n\geq 4$. He uses purely topological results and the existence of
suitable plurisubharmonic Morse functions. His method does not work
in dimensions 2 and 3. \\

\noindent {\bf Remark 5.} {\it Not every wild Cantor set in the sphere
has non-trivial polynomial hull.}\\

Indeed, this can be seen even by varying the above construction. Let the 
general step $N$ of the present construction look as follows. Suppose
we obtained sets $E^*_1\supset
\ldots \supset E^*_{N-1}$, so that the $E^*_{k}$ are disjoint unions
of closed solid tori of the same kind as above of width $2\ep^*_k$. Suppose the
polynomial hull of $E^*_{k}$ is contained in the
$\delta_k$-neighbourhood of $(\cup \ell ^{(k)})\cap
{\overline{\mathbb{B}^2}}$ for some small enough positive number
$\delta_k$. Here $\ell ^{(k)}$ denotes the collection of complex lines
which are the symmetry axes of the tori in $E^*_{k}$.
Choose at step $N$ disjoint tori $T^{*(N)} \subset E^*_{N-1} $ by
doing first the construction of the proof of proposition 4, then
fixing the symmetry axes of the tori and taking $\ep^*_N$ so
small that the polynomial hull of $E^*_N = \cup T^{*(N)}(2\ep ^*_N)$
is contained in the $\delta_N$-neighbourhood of $(\cup \ell
^{(N)})\cap {\overline{\mathbb{B}^2}}$. This is possible
since $(\cup \ell ^{(k)})\cap{\overline{\mathbb{B}^2}}$ is the
polynomial hull of $(\cup \ell ^{(k)})\cap{\partial{\mathbb{B}^2}}$.
 If $\delta_N \to 0$ for $N \to \infty $ fast enough then (since  $(\cup \ell^{(N)})\cap{\overline{\mathbb{B}^2}}$
is close to the sphere for large
$N$) the accumulation points of the mentioned
$\delta_N$-neighbourhoods are contained in the sphere. Hence the set
$E^* = \cap E^*_N$ is polynomially convex. Using Bing's theorem one
can prove that this
set is a wild Cantor set.\\

\noindent {\bf Remark 6.} {\it The set $E$ constructed in the theorem
is rationally convex.}\\

In fact, each pair of points of $E$ can be separated by a 2-torus $T$
in $\partial \mathbb{B}^2 \setminus E$. Indeed for the 2-torus $T$ we
can take the boundary of certain solid torus $T^{(N)}_m(2\ep)$, where 
$T^{(N)}_m(2\ep_N)$ is one of the tori contributing to $E_N$ for some
$N$ and $\ep$ is slightly bigger than $\ep_N$. The 2-torus $T$ is
contained in the cylinder  $\{z\in \mathbb{C}^2: |f^{(N)}_m(z)| =
2\ep\}$ which is the union of complex lines and does not meet
$E$. Hence, the cylinder does not meet the rational hull of $E$. Now
the same argument works as for polynomial convexity of tame Cantor
sets in the sphere.

We do not know whether Cantor sets in the sphere are always
rationally convex.

\nocite{*} 
\bibliography{cantor}
\bibliographystyle{plain}

\end{document}